%% file: SCOversetSRSPKMR.tex
\documentclass[a4paper,10pt]{article}

\usepackage{setspace}
\usepackage[bottom]{footmisc}
\usepackage{color}
\usepackage{subfig}
\usepackage{hyperref}
\usepackage{multirow}
\usepackage{cite}
\usepackage{amsmath} 
\usepackage{mathtools}
\usepackage{cleveref}
\usepackage{graphicx}
\usepackage{epsfig}
\usepackage{epstopdf}
\usepackage{float}
\usepackage{titling}
\usepackage[margin=1.0in]{geometry}
\usepackage{etoolbox}
 
\title{A discontinuous Galerkin overset scheme using WENO reconstruction and subcells for two-dimensional problems}
\makeatletter
\renewcommand\@date{{%
  \vspace{-\baselineskip}%
  \large\centering
  \begin{tabular}{@{}c@{}}
    S R Siva Prasad Kochi\textsuperscript{1} \\
    \normalsize siva.ksr@gmail.com
  \end{tabular}%
  \quad and\quad
  \begin{tabular}{@{}c@{}}
    M Ramakrishna\textsuperscript{2} \\
    \normalsize krishna@ae.iitm.ac.in
  \end{tabular}

  \bigskip

  \textsuperscript{1}Doctoral Candidate, Dept. of Aerospace Engg., IIT Madras.\par
  \textsuperscript{2}Professor, Dept. of Aerospace Engg., IIT Madras.

  \bigskip

  \today
}}
\makeatother
\begin{document}

\maketitle

\begin{abstract}
 A new scheme for communication between overset grids using subcells and Weighted Essentially Non Oscillatory (WENO) reconstruction for two-dimensional problems has been proposed. The effectiveness of this procedure is demonstrated using the discontinuous Galerkin method (DGM). This scheme uses WENO reconstruction using cell averages by dividing the immediate neighbors into subcells to find the degrees of freedom in cells near the overset interface. This also has the added advantage that it also works as a limiter if a discontinuity passes through the overset interface. Accuracy tests to demonstrate the maintenance of higher order are provided. Results containing shocks are also provided to demonstrate the limiter aspect of the data communication procedure.
 
 {{\bf Keywords:} overset grids, WENO reconstruction, discontinuous Galerkin method}
\end{abstract}

\section{Introduction}\label{sec:intro}

We propose a scheme for communication between overset grids using subcells and Weighted Essentially Non Oscillatory (WENO) reconstruction for two-dimensional problems. We demonstrate the effectiveness of this procedure using discontinuous Galerkin method (DGM). This scheme uses WENO reconstruction using cell averages by dividing the immediate neighbors into subcells as proposed in \cite{srspkmr1} to find the degrees of freedom in cells near the overset interface.
\\
\\
\noindent Overset meshes have been used to handle complex geometries for a long time and were first applied to solving the Euler equations by Benek et al.\cite{bsd}. A major advantage of overset meshes over a single mesh is their effectiveness for moving body problems. Overset grids communicate through exchange of boundary data (called artificial boundaries) in overlapping regions. The arbitrary overlapping of grids allows the mesh generator to focus on resolving individual components of the geometry independently. 
\\
\\
\noindent Typically Cartesian meshes are used for overset grids along with an unstructured grid in the application of high order schemes to complex geometries. In using higher order schemes, flow parameters on the artificial boundaries need to be determined so as to maintain the order of accuracy of the scheme. However, schemes like WENO \cite{shu1} require large stencils which causes problems in using artificial boundaries \cite{ss}. Discontinuous Galerkin method \cite{cs6} is well suited for overset grids as the DG discretization depends only on the current cell and its immediate neighbors. Availability of the solution polynomial in all cells is also another advantage of the DGM. DGM has been used in an overset framework for the solution of many problems in \cite{gbot}, \cite{zl}, \cite{mpfl}, \cite{bsm}.
\\
\\
\noindent When using DGM on overset grids, there are two possible approaches to handle data communication. One is a face based communication approach developed in \cite{gbot}, where solutions at an overset interface are obtained from the donor element, and then the boundary condition is applied weakly by imposing a numerical flux at the flux interpolation points. The other is an element based communication approach developed in \cite{nsm}, where the internal degrees of freedom of cells near the overset interface are obtained from the donor element.
\\
\\
\noindent In this paper, we describe a new scheme for higher order data interpolation between overset grids using the element based approach. For reconstructing the degrees of freedom in a given cell at the overset interface, we use an appropriate higher order WENO reconstruction using cell averages after dividing the immediate neighbors into subcells as proposed in \cite{srspkmr1}. This procedure has the added advantage that it also works as a limiter if a discontinuity passes through the overset interface. We can use this procedure for data communication between overset grids with any other higher order method which uses cells for their solution (eg., finite volume WENO method \cite{shu1}). We demonstrate the scheme using discontinuous Galerkin method.
\\
\\
\noindent The paper is organized as follows. We describe the formulation of the discontinuous Galerkin method used for all our results in \cref{sec:formulation}, the proposed scheme for data communication between overset grids is described in \cref{sec:oversetExpl}, and the accuracy tests and results obtained using the overset grid solver are described in \cref{sec:results} and we conclude the paper in \cref{sec:conc}.

\section{Formulation of discontinuous Galerkin Method}\label{sec:formulation}

\noindent Consider the Euler equations in conservative form as given by

\begin{equation}\label{2dEulerEquations}
\frac{\partial \textbf{Q}}{\partial t} + \frac{\partial \textbf{F(Q)}}{\partial x} + \frac{\partial \textbf{G(Q)}}{\partial y} = 0 \quad \text{in the domain} \quad \Omega
\end{equation}
\noindent where $\textbf{Q} = (\rho, \rho u, \rho v, E)^{T}$, $\textbf{F(Q)}=u\textbf{Q} + (0, p, 0, pu)^{T}$ and $\textbf{G(Q)}=v\textbf{Q} + (0, 0, p, pv)^{T}$ with $p = (\gamma -1)(E-\frac{1}{2}\rho (u^{2}+v^{2}))$ and $\gamma = 1.4$. Here, $\rho$ is the density, $(u,v)$ is the velocity, $E$ is the total energy and $p$ is the pressure. We approximate the domain $\Omega$ by $K$ non overlapping elements given by $\Omega_{k}$. 
\\
\\
We look at solving \eqref{2dEulerEquations} using the discontinuous Galerkin method. We approximate the local solution as a polynomial of order $N$ which is given by:

\begin{equation}\label{modalApprox}
 Q_{h}^{k}(r,s) = \sum_{i=0}^{N_{p}-1} Q_{i}^{k} \psi_{i}(r,s)
\end{equation}

\noindent where $N_{p}=(N+1)(N+1)$ and $r$ and $s$ are the local coordinates. The polynomial basis used ($\psi_{i}(r,s)$) is the tensor product orthonormalized Legendre polynomials of order $N$. The number of degrees of freedom are given by $N_{p}=(N+1)(N+1)$. Now, using $\psi_{j}(r,s)$ as the test function, the weak form of the equation \eqref{2dEulerEquations} is obtained as

\begin{equation}\label{weakFormScheme}
 \sum_{i=0}^{N_{p}-1} \frac{\partial Q_{i}^{k}}{\partial t} \int_{\Omega_{k}} \psi_{i} \psi_{j} d\Omega + \int_{\partial \Omega_{k}} \hat{F} \psi_{j} ds - \int_{\Omega_{k}} \vec{F}.\nabla \psi_{j} d\Omega = 0
\end{equation}

\noindent where $\partial \Omega_{k}$ is the boundary of $\Omega_{k}$, $\vec{F} = (\textbf{F(Q)},\textbf{G(Q)})$ and $\hat{F} = \bar{F^{*}}.\hat{n}$ where $\bar{F^{*}}$ is the monotone numerical flux at the interface which is calculated using an exact or approximate Riemann solver and $\hat{n}$ is the unit outward normal. This is termed to be $\mathbf{P}^{N}$ based discontinuous Galerkin method.
\\
\\
\noindent Equation \eqref{weakFormScheme} is integrated using an appropriate Gauss Legendre quadrature and is discretized in time by using the fourth order Runge-Kutta time discretization given in \cite{butcher} unless otherwise specified. To control spurious oscillations which occur near discontinuities, a limiter is used with a troubled cell indicator. We have used the KXRCF troubled cell indicator \cite{kxrcf} and the compact subcell WENO (CSWENO) limiter proposed in \cite{srspkmr1} for all our calculations.

\section{New scheme for data communication between overset grids}\label{sec:oversetExpl}

\begin{figure}[htbp]
\begin{center}
\scalebox{1.0}{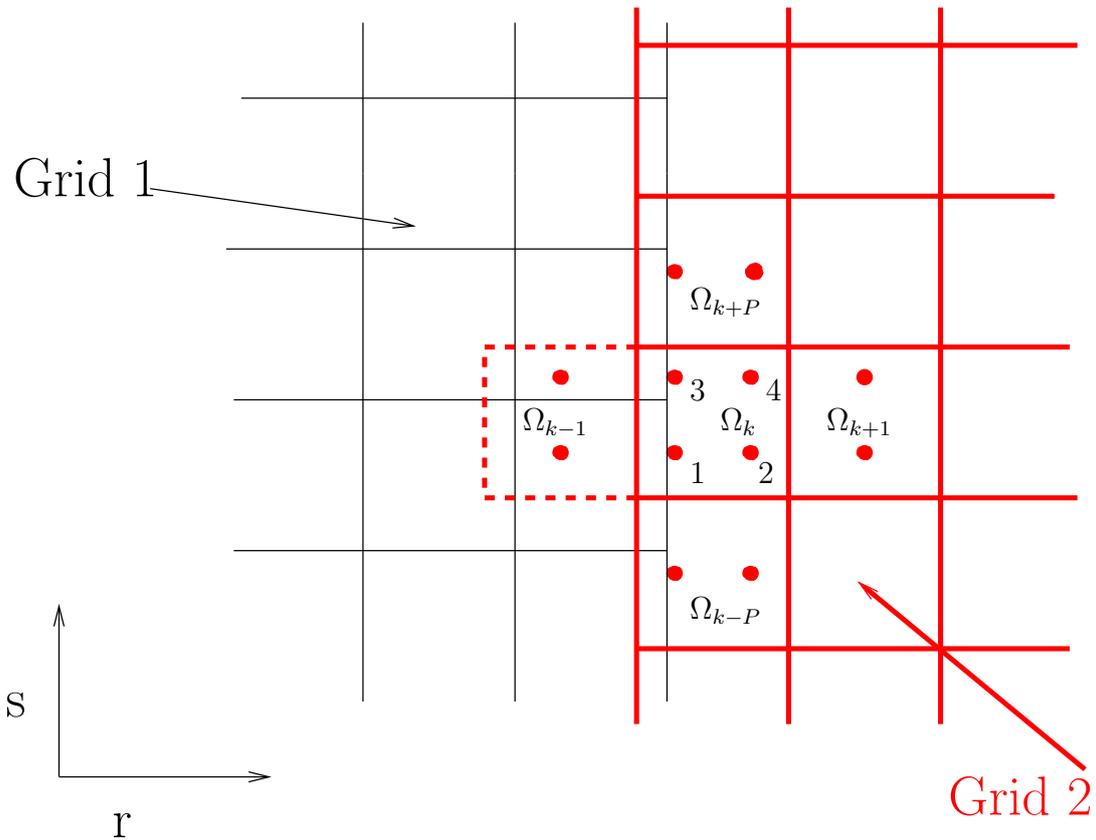}
\caption{Two overlapping grids (Grid 1 in black and Grid 2 in Red) with the element $\Omega_{k}$ where we are applying the artificial boundary conditions; $\Omega_{k}$ contains the Gauss quadrature points (points 1, 2, 3 and 4) where we reconstruct the solution for $\mathbf{P}^{1}$ based DGM; Figure also contains the ghost cell constructed for the application of artificial boundary condition (the cell $\Omega_{k-1}$) and the immediate neighbors of $\Omega_{k}$ in Grid 2 ($\Omega_{k+1}$, $\Omega_{k+P}$, $\Omega_{k-P}$) and the Gauss quadrature points in those cells used for WENO reconstruction}
\label{fig:OversetGridExpl}
\end{center}
\end{figure}

\noindent As the name suggests, overset grids consist of multiple grids which overlap each other. Boundaries of two overlapping grids named Grid 1 (black) and Grid 2 (red) are shown in Figure \ref{fig:OversetGridExpl}. These boundaries and boundary cells have an adjective artificial attached to them. For example, $\Omega_{k}$, $\Omega_{k+P}$, $\Omega_{k-P}$ are artificial boundary cells as shown in Figure \ref{fig:OversetGridExpl}. For element based data communication approach, inter-grid communication happens through the artificial boundary cells.
\\
\\

\noindent For the artificial boundary cells, we have to reconstruct the values of the degrees of freedom ($Q_{i}^{k}$ in equation \eqref{modalApprox}. In a given cell $\Omega_{k}$ in Grid 2, we keep the cell average $Q_{0}^{k}$ constant and obtain the other degrees of freedom by WENO reconstruction. To do this, we need all the immediate neighbors of $\Omega_{k}$. Across artificial boundary face of $\Omega_{k}$, we construct a ghost cell $\Omega_{k-1}$ as shown in Figure \ref{fig:OversetGridExpl} to ensure the requisite neighbors of $\Omega_{k}$ for WENO reconstruction. The ghost cell is constructed by extending the mesh line beyond the artificial boundary with the same dimension as the boundary face. This procedure for reconstruction of degrees of freedom is different from the standard element based data communication approach as it uses just the cell $\Omega_{k}$ for the reconstruction of the degrees of freedom $Q_{i}^{k}$.
\\
\\
\noindent To obtain the degrees of freedom $Q_{i}^{k}$ (i=1,2,$\ldots$,$N_{p}-1$), we use the compact subcell WENO reconstruction proposed in \cite{srspkmr1}. For $\mathbf{P^{1}}$ based DGM, we reconstruct the values of $Q$ at points $1,2,3,4$ shown in Figure \ref{fig:OversetGridExpl}. These points correspond to the tensor product of the two Gauss quadrature points. We first calculate $Q^{k-1}_{h}$ and $Q^{k+1}_{h}$ at coordinates $(r=0,s=\pm 1/\sqrt{3})$ and $Q^{k-P}_{h}$ and $Q^{k+P}_{h}$ at coordinates $(r=\pm 1/\sqrt{3},s=0)$ corresponding to the points shown in Figure \ref{fig:OversetGridExpl}.
\\
\\
\noindent In this example, all the required values in cells $\Omega_{k+1}$, $\Omega_{k+P}$, and $\Omega_{k-P}$ are calculated using the cell local DG polynomial in Grid 2 using equation \eqref{modalApprox}. We explain how to find the values of $Q$ in the ghost cell $\Omega_{k-1}$ now. Assume that the points where we need the values of $Q$ (the points $r=0,s=\pm 1/\sqrt{3}$) are represented in Cartesian coordinates as $(x_{1}(r_{1},s_{1}),y_{1}(r_{1},s_{1}))$ and $(x_{2}(r_{2},s_{2}),y_{2}(r_{2},s_{2}))$, where $(r_{1},s_{1})$ and $r_{2},s_{2})$ are the point locations in local coordinate system. These are used to find the Cartesian coordinates $\bar{X_{k}}$ for $k=1,2$. These Cartesian coordinates are used to obtain the cell local coordinates in Grid 1 which are $(r_{k}(\bar{X_{k}}),s_{k}(\bar{X_{k}}))$ for $k=1,2$. A search algorithm (K-d tree) \cite{dcv} is used to determine which cell of Grid 1 contains the given GQ points. Then, we can find the cell local coordinates in a given cell $\Omega_{i}$ using Newton's method as given below:

\begin{equation}\label{NewtonMethod}
 \left[\begin{array}{cc}\partial_{r}x_{i}(r^{n},s^{n}) & \partial_{s}x_{i}(r^{n},s^{n}) \\ \partial_{r}y_{i}(r^{n},s^{n}) & \partial_{s}y_{i}(r^{n},s^{n})\end{array}\right]\left[\begin{array}{cc}\Delta r \\ \Delta s\end{array}\right]=\left[\begin{array}{cc} -(x_{i}(r^{n},s^{n})-x_{b}(1,s_{k})) \\ -(y_{i}(r^{n},s^{n})-y_{b}(1,s_{k}))\end{array}\right] \quad \text{for  } k = 1,2
\end{equation}

\noindent where $r^{0}=0$, $s^{0}=0$, $r^{n+1}=r^{n}+\Delta r$, $s^{n+1}=s^{n}+\Delta s$, $\partial_{r}=\partial/\partial r$, and $\partial_{s}=\partial/\partial s$. We iterate using \eqref{NewtonMethod} until the Euclidean distance between the terms on the right hand side of \eqref{NewtonMethod} drops below a tolerance of $1e-10$ or the Newton method reaches a maximum number of 20 iterations. The cell $\Omega_{i}$ is a donor cell for the coordinate $\bar{X_{k}}$ if the Euclidean distance drops below the tolerance of $1e-10$. If the Euclidean distance is above the required tolerance after 20 iterations, the coordinate $\bar{X_{k}}$ is deemed to reside outside of the cell $\Omega_{i}$. In this way, we obtain the cell local coordinates in Grid 1 corresponding to $(x_{1}(r_{1},s_{1}),y_{1}(r_{1},s_{1}))$ and $(x_{2}(r_{2},s_{2}),y_{2}(r_{2},s_{2}))$.
\\
\\
\noindent If the mesh changes to adapt to the solution or the changing geometry, these coordinates have to be calculated at every time step. Otherwise, they are calculated and stored before hand so as not to calculate them at every time step for stationary overset grids. Using these cell local coordinates and the DG polynomial in that cell, we can find the $Q$ values at each of the required points using equation \eqref{modalApprox}. 
\\
\\
\noindent For our reconstruction, we need the values of $Q$ at locations $(r=0,s=\pm 1/\sqrt{3})$ in cells $\Omega_{k-1}$ and $\Omega_{k+1}$ and at $(r=\pm 1/\sqrt{3},s=0)$ in cells $\Omega_{k-P}$ and $\Omega_{k+P}$. After finding the values of $Q$ at these points which are shown in Figure \ref{fig:OversetGridExpl} in each cell, we use $Q^{k-1}_{h}(r=0,s=-1/\sqrt{3})$, $Q_{0}^{k}$ and $Q^{k+1}_{h}(r=0,s=-1/\sqrt{3})$ to find $Q_{x1}$ and $Q_{x2}$ with a WENO3 reconstruction. Similarly, $Q^{k-1}_{h}(r=0,s=1/\sqrt{3})$, $Q_{0}^{k}$ and $Q^{k+1}_{h}(r=0,s=1/\sqrt{3})$ are used to find $Q_{x3}$ and $Q_{x4}$. Here, $Q_{x1}$, $Q_{x2}$, $Q_{x3}$ and $Q_{x4}$ are the one-dimensional WENO reconstructed values at the Gauss quadrature points $1$, $2$, $3$ and $4$ shown in Figure \ref{fig:OversetGridExpl} in the $x$ direction. In the same manner, we use $Q^{k-P}_{h}(r=\pm 1/\sqrt{3},s=0)$, $Q^{k+P}_{h}(r=\pm 1/\sqrt{3},s=0)$, and $Q_{0}^{k}$ to find $Q_{y1}$, $Q_{y2}$, $Q_{y3}$ and $Q_{y4}$. Again, $Q_{y1}$, $Q_{y2}$, $Q_{y3}$ and $Q_{y4}$ are the one-dimensional WENO reconstructed values at the Gauss quadrature points $1$, $2$, $3$ and $4$ in the $y$ direction.
\\
\\
\noindent We now describe the procedure to obtain $Q_{x1}$, $Q_{x2}$, $Q_{x3}$, $Q_{x4}$, $Q_{y1}$, $Q_{y2}$, $Q_{y3}$ and $Q_{y4}$ using one-dimensional WENO3 reconstruction briefly. This is described in detail in \cite{srspkmr1}. We describe the procedure to find the third order WENO reconstruction to obtain the values of $Q$ in cell $\Omega_{k}$ at the Gauss quadrature points $Q_{x1}$ and $Q_{x2}$, given three cells $\Omega_{k-1}$, $\Omega_{k}$, and $\Omega_{k+1}$, and the corresponding cell averages $Q_{0}^{k-1}$, $Q_{0}^{k}$ and $Q_{0}^{k+1}$. We identify $2$ small stencils $S_{i}$, $i=0,1$ such that $\Omega_{k}$ belongs to each of them. We set $S_{i} = \bigcup_{l=0}^{1}\Omega_{k+i-l}$. We also have the larger stencil $\mathbf{T}=\bigcup_{i=0}^{1}S_{i}$ which contains all the cells from the smaller stencils $S_{i}$.
\\
\\
\noindent Now, we have a polynomial of degree $1$, $p_{i}(x)$ corresponding to the stencil $S_{i}$ such that it's cell average in each of the cells of the stencil $S_{i}$ agrees with the given cell average of $Q$. We also have a polynomial of degree $2N$ reconstruction denoted by $M(x)$ associated with the larger stencil $\mathbf{T}$, such that the cell average of $M(x)$ in each of the cells of the stencil $\mathbf{T}$ agrees with the cell average of $Q$ for that cell. The details of the construction of $p_{i}(x)$ and $M(x)$ are given in \cite{shu1}.
\\
\\
\noindent Next, we find the linear weights denoted by $\gamma_{0},\gamma_{1}$, which satisfy

\begin{equation}\label{Chap3Eq2:WENO1DReconstructionLinearWeightsDefinition}
 M(x_{G}) = \sum_{i=0}^{1} \gamma_{i} p_{i}(x_{G})
\end{equation}

\noindent where $x_{G}$ is a Gauss or Gauss-Lobatto quadrature point. A set of linear weights for each of the quadrature points is obtained. The value of the functions $M(x)$ and $p_{i}(x)$ for each $i$ can be written as a function of the cell average of each cell in the stencil. This is used in WENO reconstruction. For the $P^{1}$ based DGM, with the Gauss quadrature point $r=-1/\sqrt{3}$, we have:

\begin{equation}\label{Chap3Eq3:WENO1DReconstructionPolynomialsP1_1}
 \left[L^{1}_{i}\right] = \left[T^{1}_{ij}\right] \left[C^{1}_{j}\right]
\end{equation}

\noindent where

\begin{displaymath}
\left[L^{1}_{i}\right] =  \begin{bmatrix} p_{0}(x_{G}) & p_{1}(x_{G}) & M(x_{G}) \end{bmatrix}^{T}
\end{displaymath}

\begin{displaymath}
 \left[C^{1}_{j}\right] = \begin{bmatrix} Q_{0}^{k-1} & Q_{0}^{k} & Q_{0}^{k+1} \end{bmatrix}^{T}
\end{displaymath}

\noindent and

\begin{displaymath}
 \left[T^{1}_{ij}\right] = \begin{bmatrix} \frac{\sqrt{3}}{6} & \frac{6-\sqrt{3}}{6} & 0 \\ 0 & \frac{6+\sqrt{3}}{6} & -\frac{\sqrt{3}}{6} \\ \rule{0pt}{2.5ex}  \frac{\sqrt{3}}{12} & 1 & \rule{0pt}{2.5ex} -\frac{\sqrt{3}}{12} \end{bmatrix}
\end{displaymath}

\noindent The linear weights are given by

\begin{equation}\label{Chap3Eq4:WENO1DReconstructionLinearWeightsP1_1}
 \gamma_{0} = \frac{1}{2}, \qquad \qquad \gamma_{1} = \frac{1}{2}
\end{equation}

\noindent For the Gauss quadrature point $r=1/\sqrt{3}$, we have:

\begin{equation}\label{Chap3Eq3:WENO1DReconstructionPolynomialsP1_2}
 \left[L^{1}_{i}\right] = \left[T^{2}_{ij}\right] \left[C^{1}_{j}\right]
\end{equation}

\noindent where 

\begin{displaymath}
 \left[T^{2}_{ij}\right] = \begin{bmatrix} -\frac{\sqrt{3}}{6} & \frac{6+\sqrt{3}}{6} & 0 \\ \rule{0pt}{2.5ex} 0 & \rule{0pt}{2.5ex} \frac{6-\sqrt{3}}{6} & \rule{0pt}{2.5ex} \frac{\sqrt{3}}{6} \\ \rule{0pt}{2.5ex} -\frac{\sqrt{3}}{12} & \rule{0pt}{2.5ex} 1 & \rule{0pt}{2.5ex} \frac{\sqrt{3}}{12} \end{bmatrix}
\end{displaymath}

\noindent The linear weights remain the same.
\\
\\
\noindent Now, as given by \cite{js1}, we compute the smoothness indicator for each stencil $S_{i}$:

\begin{equation}\label{Chap3Eq4:WENO1DReconstructionSmoothnessIndicatorDefinition}
\beta_{i} = \sum_{l=1}^{N}\int_{\Omega_{j}} \Delta x_{j}^{2l-1} \left(\frac{\partial^{l}}{\partial x^{l}}p_{i}(x)\right)^{2} dx
\end{equation}

\noindent For $P^{1}$ based DGM, the smoothness indicators are given as:

\begin{equation}\label{Chap3Eq4:WENO1DReconstructionSmoothnessIndicatorP1_1}
 \beta_{0} = (Q_{0}^{k}-Q_{0}^{k-1})^{2}
\end{equation}

\begin{equation}\label{Chap3Eq4:WENO1DReconstructionSmoothnessIndicatorP1_2}
 \beta_{1} = (Q_{0}^{k+1}-Q_{0}^{k})^{2}
\end{equation}

\noindent Now, we compute the nonlinear weights as given below:

\begin{equation}\label{Chap3Eq4:WENO1DReconstructionNonLinearWeightsDefinition}
\omega_{i}=\frac{\bar{\omega}_{i}}{\sum_{i}\bar{\omega}_{i}}, \quad \bar{\omega_{i}} = \frac{\gamma_{i}}{\sum_{i}(\epsilon + \beta_{i})^{2}}
\end{equation}

\noindent Here $\epsilon$ is a small number which is usually taken to be $10^{-6}$. The final WENO approximation is given by

\begin{equation}\label{Chap3Eq4:WENO1DReconstructionWENOApprox}
Q_{G} = \sum_{i=0}^{N}\omega_{i}p_{i}(x_{G})
\end{equation}

\begin{figure}[htbp]
\begin{center}
\scalebox{1.0}{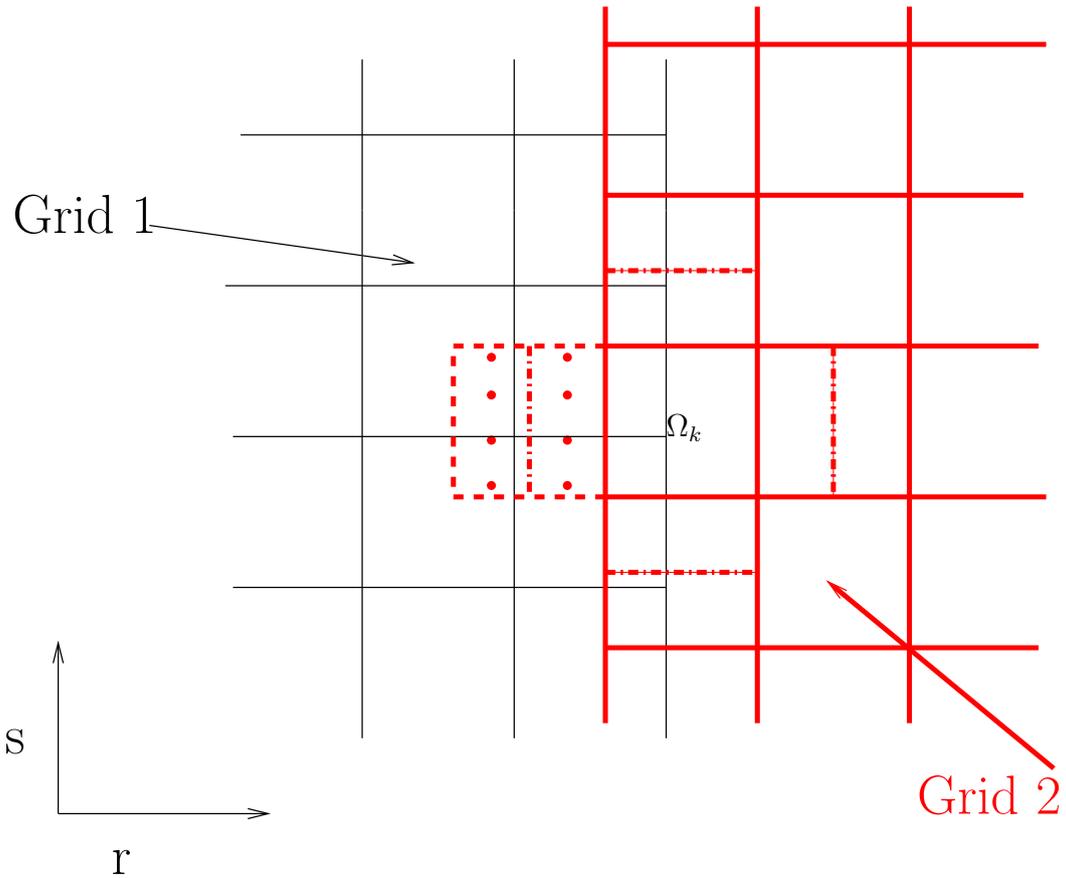}
\caption{Two overlapping grids (Grid 1 in black and Grid 2 in Red) with the element $\Omega_{k}$ where we are applying the artificial boundary conditions; The immediate neighbors divided in half for WENO5 reconstruction (used for $\mathbf{P}^{2}$ based DGM) are shown and the Gauss quadrature points in the ghost cell used for WENO reconstruction is also shown}
\label{fig:OversetGridExpl2}
\end{center}
\end{figure}

\noindent Finally, we obtain the reconstructed degrees of freedom based on the reconstructed point values $Q(x_{G})$ at the Gauss quadrature points $x_{G}$ and a numerical integration as

\begin{equation}\label{Chap3Eq4:WENO1DReconstructionNumericalIntegration}
Q_{j}^{k} = \Delta x_{k} \sum_{G} w_{G} Q(x_{G}) \psi_{j}(x_{G}) \quad j=1,\ldots,N
\end{equation}

\noindent where $w_{G}$'s are the Gaussian quadrature weights for the points $x_{G}$. This procedure allows us to find all the values $Q_{x1}$, $Q_{x2}$, $Q_{x3}$, $Q_{x4}$, $Q_{y1}$, $Q_{y2}$, $Q_{y3}$ and $Q_{y4}$ while maintaining the order of the scheme as described in \cite{srspkmr1}.
\\
\\
\noindent Using $Q_{x1}$, $Q_{x2}$, $Q_{x3}$ and $Q_{x4}$, we get the degrees of freedom $Q_{xn}^{k}$ corresponding to a polynomial in the $x$ direction $\forall n=1\ldots N_{p}-1$. Similarly, we get the degrees of freedom $Q_{yn}^{k}$ corresponding to a polynomial in the $y$ direction $\forall n=1\ldots N_{p}-1$. Now, we use the scheme

\begin{align}\label{weakFormScheme2DNew}
 \sum_{i=0}^{N_{p}-1} \frac{\partial Q_{i}^{kNew}}{\partial t} \int_{\Omega_{k}} \psi_{i} \psi_{j} d\Omega + \int_{\partial \Omega_{k}} \hat{F}(Q_{xi}^{k}, Q_{yi}^{k})  \psi_{j} ds - \int_{\Omega_{k}} \vec{F}(Q_{xi}^{k}, Q_{yi}^{k}).\nabla \psi_{j} d\Omega = 0
\end{align}
\noindent where the fluxes $\vec{F}$ and $\hat{F}$ are calculated using the appropriate values of $Q_{xn}^{k}$ and $Q_{yn}^{k}$. We also use ($Q_{hx}^{k}$ + $Q_{hy}^{k}$)/2 as $Q_{h}^{kNew}$ for time integration. For solving a system of equations, we use this with a local characteristic field decomposition with the corresponding Jacobians in the $x$ and $y$ directions as explained in \cite{zs}. This completes the procedure for data communication between overset grids in a given artificial boundary cell. We repeat this procedure for all boundary cells in both Grid 1 and Grid 2.
\\
\\
\noindent For $\mathbf{P}^{2}$ based DGM, we follow the same procedure after dividing the immediate neighbors $\Omega_{k-1}$, $\Omega_{k+1}$, $\Omega_{k-P}$ and $\Omega_{k+P}$ in half as shown in Figure \ref{fig:OversetGridExpl2} and assigning appropriate values to the new cells as given in \cite{srspkmr1}. For reconstruction of degrees of freedom with $\mathbf{P}^{2}$ based DGM, we have to use the four point Gauss-Lobatto quadrature as the corresponding Gauss quadrature rule requires the point $r=0$, where the reconstruction of the solution loses it's accuracy. The Gauss-Lobatto quadrature points used in the ghost cell are shown in Figure \ref{fig:OversetGridExpl2}. Now, we follow the procedure for WENO5 reconstruction as given in \cite{srspkmr1} and \cite{shu1}. For $\mathbf{P}^{3}$ based DGM, we follow the same procedure of constructing ghost cells and subcells using four point Gauss quadrature.
\\
\\

\noindent We described a scheme for data communication using WENO reconstruction and subcells for overset grids. This procedure also works as a limiter for artificial boundary cells if a discontinuity passes through the overlapping region while maintaining the order of accuracy of the solution. This makes this scheme quite effective as limiting in an artificial boundary cell is quite difficult as most of the current limiting procedures are difficult to apply on a boundary.

\section{Results}\label{sec:results}

\noindent In this section, we look at some of the results obtained to demonstrate the performance of the scheme for data communication described above. All the results are obtained using DG method and an appropriate Runge-Kutta scheme for time integration \cite{butcher}.

\subsection{Accuracy Tests}\label{subsec:accuracyTest}

\begin{figure}[htbp]
\begin{center}
\includegraphics[scale=0.4]{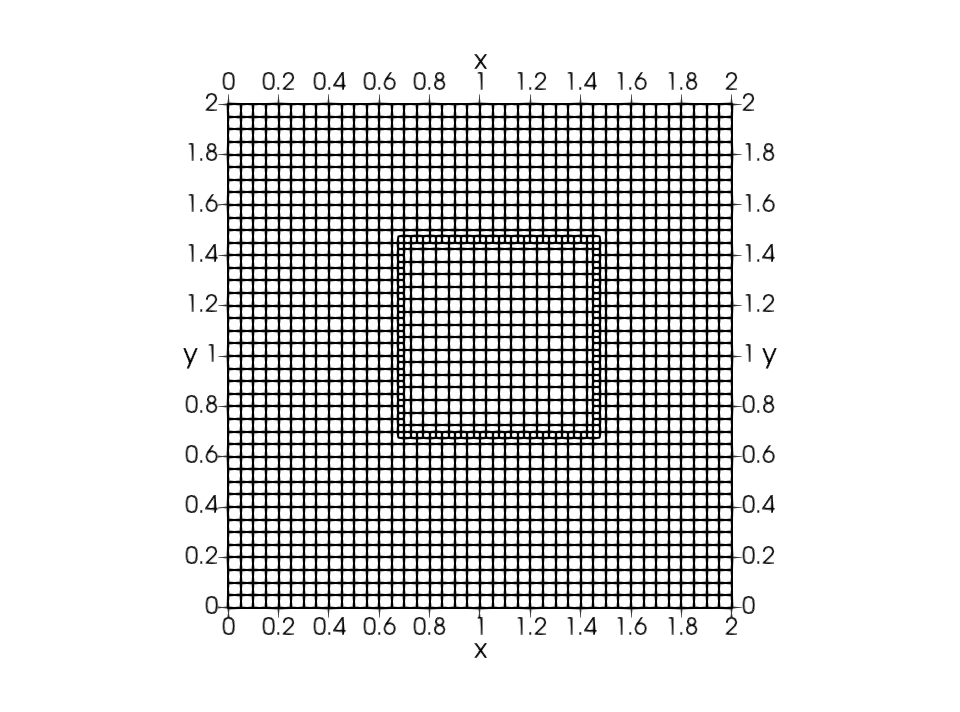}
\caption{Overset Grid of Type 1 for Entropy Wave Problem in the domain $[0,2]\times[0,2]$ with $40$ by $40$ elements used for the validation of the overset grid solver.}
\label{fig:EntropyWaveOversetGrid1}
\end{center}
\end{figure}

\begin{figure}[htbp]
\begin{center}
\includegraphics[scale=0.4]{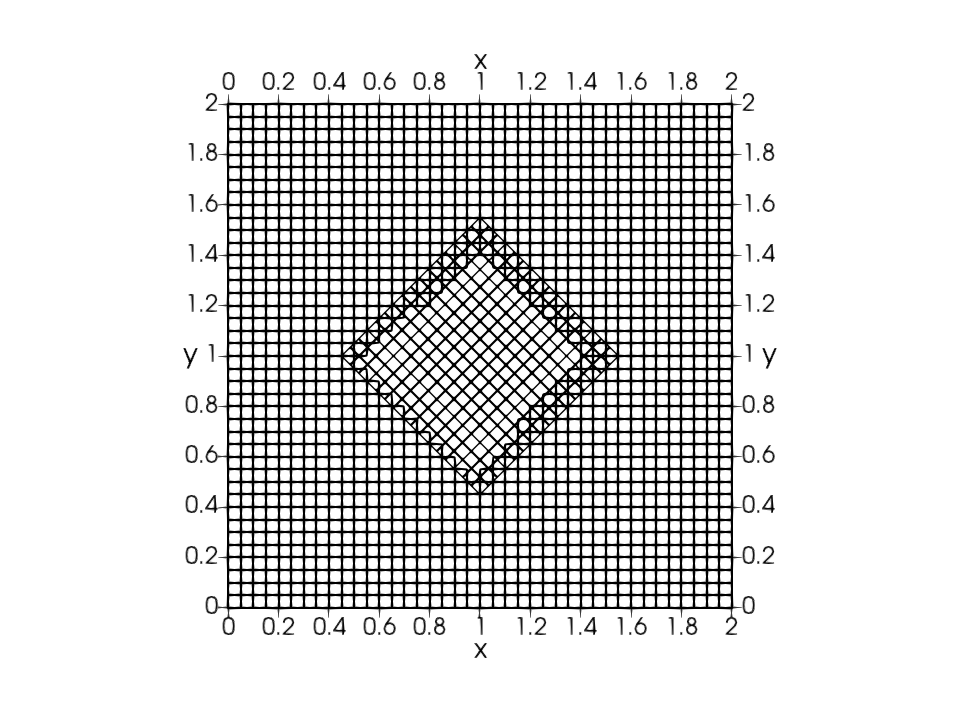}
\caption{Overset Grid of Type 2 (which is obtained by rotating the Type 1 overset grid by $45^{0}$) for Entropy Wave Problem in the domain $[0,2]\times[0,2]$ with $40$ by $40$ elements used for the validation of the overset grid solver.}
\label{fig:EntropyWaveOversetGrid2}
\end{center}
\end{figure}

\noindent \textbf{Example 1:} We solve the two dimensional Euler equations as given in equation \eqref{2dEulerEquations} in the domain $[0,2]\times[0,2]$. The initial conditions are given by $\rho(x,y,0) = 1+0.2\sin(\pi (x+y))$, $u(x,y,0) = 0.7$, $v(x,y,0) = 0.3$ and $p(x,y,0) = 1.0$ and we use periodic boundary conditions in both directions. The exact solution is given by $\rho(x,y,0) = 1+0.2\sin(\pi (x+y-t))$, $u(x,y,0) = 0.7$, $v(x,y,0) = 0.3$ and $p(x,y,0) = 1.0$. We run the solver with the same grid size for a normal grid and an overset grid for grid sizes of $1/20$, $1/40$, $1/80$, and $1/160$ for various orders. We have used two different overset grids as shown in Figures \ref{fig:EntropyWaveOversetGrid1} (labeled as Type 1) and \ref{fig:EntropyWaveOversetGrid2} (labeled as Type 2 which is obtained by rotating the Type 1 overset grid by $45^{0}$) for our calculations to demonstrate the effectiveness of our procedure. Both example grids shown in Figures \ref{fig:EntropyWaveOversetGrid1} and \ref{fig:EntropyWaveOversetGrid2} contain $40$ by $40$ elements. The errors in density and numerical orders of accuracy are calculated at $t=2.0$ for the original grid as well as both the overset grids and are presented in Table \ref{table:1}. While calculating the solution, we have made sure that the temporal and spatial orders of accuracy are the same by using a corresponding Runge-Kutta time integration \cite{butcher}. We can see that the solution obtained using the overset grid is as accurate as the solution obtained without any overset.

\begin{table}
\small
\centering
\resizebox{\textwidth}{!}{%
\begin{tabular}{|c|c|c|c|c|c|c|c|}
\hline
\multirow{2}{*}{} &  & \multicolumn{2}{|c|}{DG w/o overset} & \multicolumn{2}{|c|}{DG with overset (Type 1)} & \multicolumn{2}{|c|}{DG with overset (Type 2)}\\ \cline{2-8} 
 & Grid size & $L_{2}$ error & Order & $L_{2}$ error & Order & $L_{2}$ error & Order \\ \hline
 \multirow{4}{*}{$\mathbf{P}^{1}$} & 1/20 & 1.384E-03 &  & 9.976E-04 &  & 1.174E-03 & \\ \cline{2-8}
  & 1/40 & 2.850E-04 & 2.28 & 2.026E-04 & 2.30 & 2.401E-04 & 2.29 \\ \cline{2-8}
  & 1/80 & 6.289E-05 & 2.18 & 4.349E-05 & 2.22 & 5.262E-05 & 2.19 \\ \cline{2-8}
  & 1/160 & 1.437E-05 & 2.13 & 9.867E-06 & 2.14 & 1.202E-05 & 2.13 \\ \hline
  \multirow{4}{*}{$\mathbf{P}^{2}$} & 1/20 & 1.087E-05 &  & 9.876E-06 &  & 1.021E-05 &  \\ \cline{2-8}
  & 1/40 & 1.127E-06 & 3.27 & 9.889E-07 & 3.32 & 1.051E-06 & 3.28 \\ \cline{2-8}
  & 1/80 & 1.185E-07 & 3.25 & 1.018E-07 & 3.28 & 1.097E-07 & 3.26 \\ \cline{2-8}
  & 1/160 & 1.272E-08 & 3.22 & 1.085E-08 & 3.23 & 1.177E-08 & 3.22 \\ \hline
  \multirow{4}{*}{$\mathbf{P}^{3}$} & 1/20 & 1.004E-07 &  & 9.157E-08 &  & 9.938E-08 & \\ \cline{2-8}
  & 1/40 & 4.690E-09 & 4.42 & 4.337E-09 & 4.40 & 4.675E-09 & 4.41 \\ \cline{2-8}
  & 1/80 & 2.252E-10 & 4.38 & 2.069E-10 & 4.39 & 2.245E-10 & 4.38 \\ \cline{2-8}
  & 1/160 & 1.135E-11 & 4.31 & 1.029E-11 & 4.33 & 1.124E-11 & 4.32 \\ \hline
\end{tabular}}
\caption{Validation of overset grid solver using 2D Euler equations for the Entropy Wave problem with periodic boundary conditions, $t=2.0$, Uniform mesh with and without overset for two different overset grids as shown in Figures \ref{fig:EntropyWaveOversetGrid1} and \ref{fig:EntropyWaveOversetGrid2}, $L_{2}$ error for density with $\mathbf{P}^{1}$, $\mathbf{P}^{2}$ and $\mathbf{P}^{3}$ based DGM}
\label{table:1}
\end{table}

\begin{figure}[htbp]
\begin{center}
\includegraphics[scale=0.4]{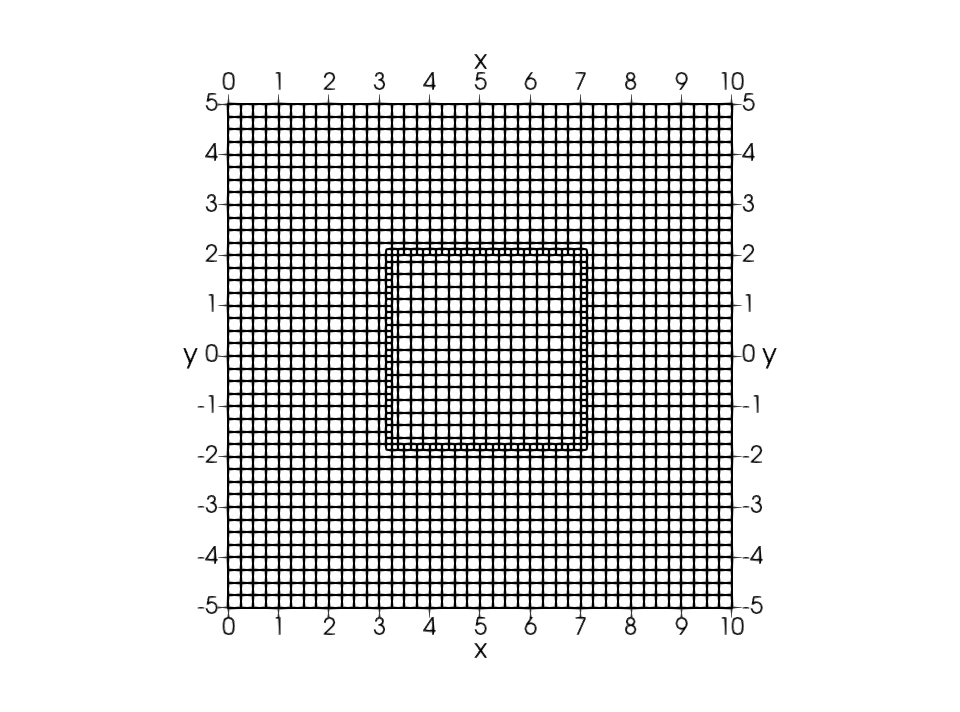}
\caption{Overset Grid of Type 1 for Isentropic Vortex Problem in the domain $[0,10]\times[-5,5]$ with $40$ by $40$ elements used for the validation of the overset grid solver.}
\label{fig:IsenVortexOversetGrid1}
\end{center}
\end{figure}

\begin{figure}[htbp]
\begin{center}
\includegraphics[scale=0.4]{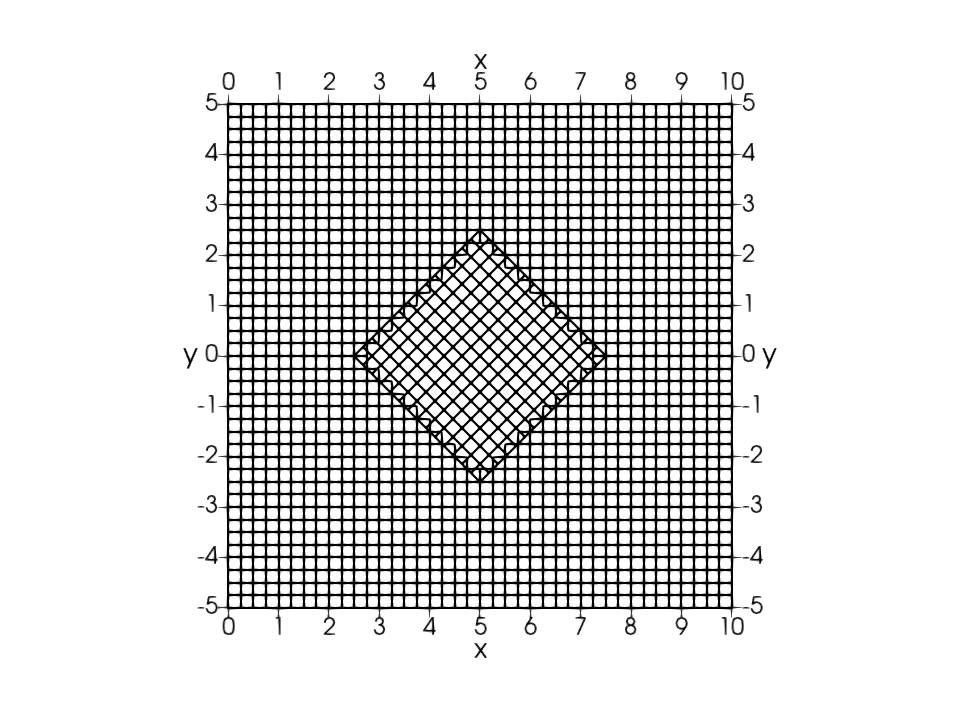}
\caption{Overset Grid of Type 2 (which is obtained by rotating the Type 1 overset grid by $45^{0}$) for Isentropic Vortex Problem in the domain $[0,10]\times[-5,5]$ with $40$ by $40$ elements used for the validation of the overset grid solver.}
\label{fig:IsenVortexOversetGrid2}
\end{center}
\end{figure}

\noindent \textbf{Example 2:} We solve the two dimensional Euler equations as given in equation \eqref{2dEulerEquations} in the domain $[0,10]\times[-5,5]$ for the Isentropic Euler Vortex problem. The analytical solution is given by: \\ $\rho = \left(1 -  \left(\frac{\gamma - 1}{16\gamma \pi^{2}}\right)\beta^{2} e^{2(1-r^{2})}\right)^{\frac{1}{\gamma-1}}$, $u = 1 - \beta e^{(1-r^{2})} \frac{y-y_{0}}{2\pi}$, $v = \beta e^{(1-r^{2})} \frac{x-x_{0}-t}{2\pi}$, and $p = \rho^{\gamma}$, where $r$ is given by $\sqrt{(x-x_{0}-t)^{2}+(y-y_{0})^{2}}$, $x_{0}=5$, $y_{0}=0$, $\beta=5$ and $\gamma = 1.4$. We initialize with the analytical solution at $t=0$ and use periodic boundary conditions at the edges of the domain. We run the solver with the same grid size for the baseline grid and an overset grid for grid sizes of $1/20$, $1/40$, $1/80$, and $1/160$ for various orders. We have used two different overset grids as shown in Figures \ref{fig:IsenVortexOversetGrid1} (labeled as Type 1) and \ref{fig:IsenVortexOversetGrid2} (labeled as Type 2 which is obtained by rotating the Type 1 overset grid by $45^{0}$) for our calculations to demonstrate the effectiveness of our procedure. Both example grids shown in Figures \ref{fig:IsenVortexOversetGrid1} and \ref{fig:IsenVortexOversetGrid2} contain $40$ by $40$ elements. The errors in density and numerical orders of accuracy are calculated at $t=10.0$ (one period) for the original grid as well as both the overset grids and are presented in Table \ref{table:2}. While calculating the solution, we have made sure that the temporal and spatial orders of accuracy are the same by using a corresponding Runge-Kutta time integration \cite{butcher}. We can see that the solution obtained using the overset grid is as accurate as the solution obtained without any overset.

\begin{table}
\small
\centering
\resizebox{\textwidth}{!}{%
\begin{tabular}{|c|c|c|c|c|c|c|c|}
\hline
\multirow{2}{*}{} &  & \multicolumn{2}{|c|}{DG w/o overset} & \multicolumn{2}{|c|}{DG with overset (Type 1)} & \multicolumn{2}{|c|}{DG with overset (Type 2)}\\ \cline{2-8} 
 & Grid size & $L_{2}$ error & Order & $L_{2}$ error & Order & $L_{2}$ error & Order \\ \hline
 \multirow{4}{*}{$\mathbf{P}^{1}$} & 1/20 & 3.215E-03 &  & 1.252E-03 &  & 2.786E-03 & \\ \cline{2-8}
  & 1/40 & 7.294E-04 & 2.14 & 2.801E-04 & 2.16 & 6.277E-04 & 2.15 \\ \cline{2-8}
  & 1/80 & 1.725E-04 & 2.08 & 6.534E-05 & 2.10 & 1.485E-04 & 2.08 \\ \cline{2-8}
  & 1/160 & 4.137E-05 & 2.06 & 1.535E-05 & 2.09 & 3.561E-05 & 2.06 \\ \hline
  \multirow{4}{*}{$\mathbf{P}^{2}$} & 1/20 & 2.232E-05 &  & 1.493E-05 &  & 1.876E-05 &  \\ \cline{2-8}
  & 1/40 & 2.512E-06 & 3.15 & 1.647E-06 & 3.18 & 2.099E-06 & 3.16 \\ \cline{2-8}
  & 1/80 & 2.829E-07 & 3.15 & 1.830E-07 & 3.17 & 2.365E-07 & 3.15 \\ \cline{2-8}
  & 1/160 & 3.187E-08 & 3.15 & 2.047E-08 & 3.16 & 2.664E-08 & 3.15 \\ \hline
  \multirow{4}{*}{$\mathbf{P}^{3}$} & 1/20 & 2.768E-07 &  & 1.842E-07 &  & 2.274E-07 & \\ \cline{2-8}
  & 1/40 & 1.376E-08 & 4.33 & 8.847E-09 & 4.38 & 1.123E-08 & 4.34 \\ \cline{2-8}
  & 1/80 & 6.747E-10 & 4.35 & 4.279E-10 & 4.37 & 5.545E-10 & 4.34 \\ \cline{2-8}
  & 1/160 & 3.378E-11 & 4.32 & 2.128E-11 & 4.33 & 2.796E-11 & 4.31 \\ \hline
\end{tabular}}
\caption{Validation of overset grid solver using 2D Euler equations for the Isentropic Vortex problem with periodic boundary conditions, $t=10.0$,Uniform mesh with and without overset for two different overset grids as shown in Figures \ref{fig:IsenVortexOversetGrid1} and \ref{fig:IsenVortexOversetGrid2}, $L^{2}$ error for density with $\mathbf{P}^{1}$, $\mathbf{P}^{2}$ and $\mathbf{P}^{3}$ based DGM}
\label{table:2}
\end{table}

\subsection{Test Cases with discontinuities}\label{subsec:ResultsWithDiscontinuities}

\begin{figure}[htbp]
\begin{center}
\includegraphics[scale=0.2]{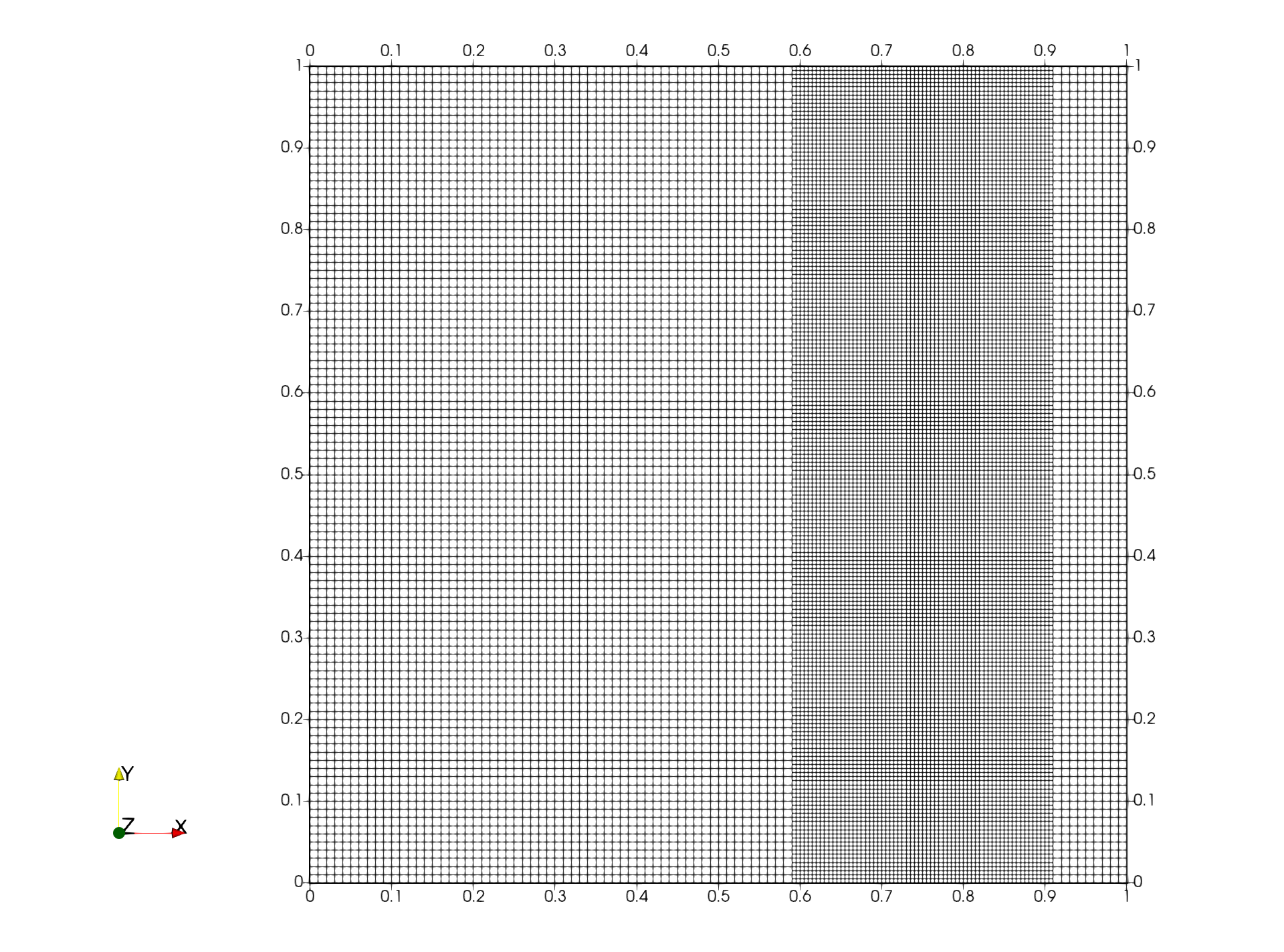}
\caption{Overset Grid used for the solution of Sod's shock tube problem in the domain $[0,1]\times[0,1]$ with a refined overset grid of size $h=1/200$ between $x=0.59$ and $x=0.91$ (where the solution contains a discontinuity) on a baseline grid of size $h=1/100$.}
\label{fig:SodOversetGrid}
\end{center}
\end{figure}

\begin{figure}[htbp]
\begin{center}
\includegraphics[scale=0.7]{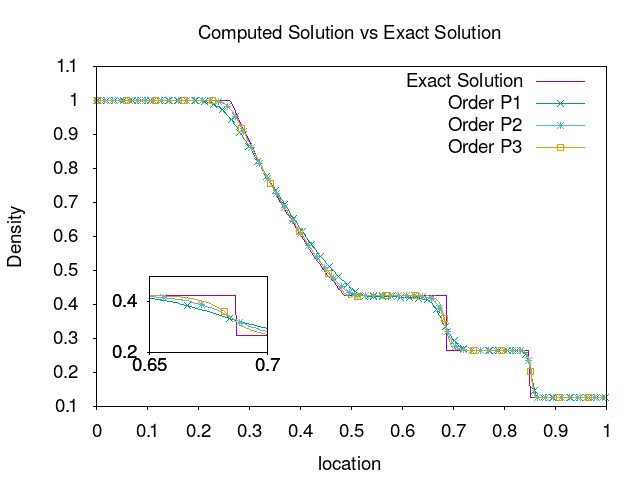}
\caption{Comparison of density solutions on $y=0.5$ line of Sod's shock tube problem at $t=0.2$ in the domain $[0,1]\times[0,1]$ with a refined overset grid of size $h=1/200$ between $x=0.59$ and $x=0.91$ (where the solution contains a discontinuity) on a baseline grid of size $h=1/100$ obtained with the $P^{1}$, $P^{2}$ and $P^{3}$ based DGM using the proposed data communication scheme and the exact solution. Figure also includes a zoomed in portion of the solution for better comparison}
\label{fig:SodOversetSolution}
\end{center}
\end{figure}

\begin{figure}[htbp]
\begin{center}
\includegraphics[scale=0.7]{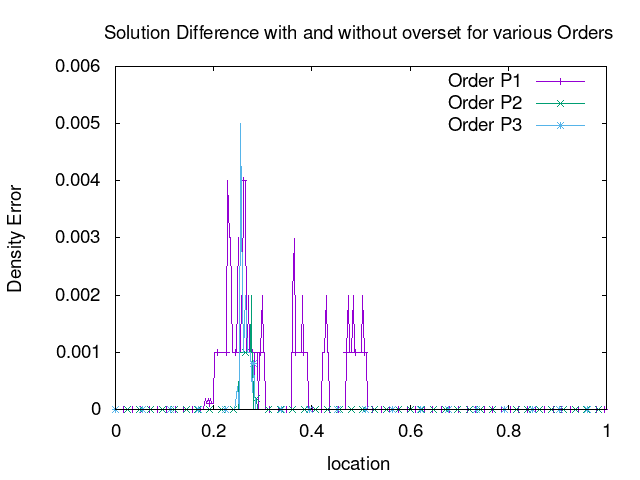}
\caption{Solution difference for density ($|\rho_{withOverset}-\rho_{withOutOverset}|$) obtained on $y=0.5$ line of Sod's shock tube problem at $t=0.2$ in the domain $[0,1]\times[0,1]$ for $P^{1}$, $P^{2}$ and $P^{3}$ based DGM between a refined overset grid of size $h=1/200$ between $x=0.59$ and $x=0.91$ (where the solution contains a discontinuity) on a baseline grid of size $h=1/100$ and a single grid of size $h=1/200$.}
\label{fig:SodOversetError}
\end{center}
\end{figure}

\noindent We now test the data communication scheme for problems with solutions having discontinuities, some of them passing through the artificial boundary. We have used the compact subcell WENO limiter proposed in \cite{srspkmr1} along with the KXRCF troubled cell indicator \cite{kxrcf} for all our calculations.
\\
\\
\noindent \textbf{Example 3:} We solve the Sod's shock tube problem as proposed in \cite{sod} in the two-dimensional domain. We solve the 2D Euler equations in the domain $[0,1] \times [0,1]$ with the initial conditions given as $(\rho,u,v,p)=(1.0,0.0,0.0,1.0)$ for $x < 0.5$ and $(\rho,u,v,p)=(0.125,0.0,0.0,0.1)$ otherwise. Non reflecting boundary condition is applied at $x=0$ and $x=1$ and periodic boundary conditions are applied at the other two boundaries. We use a refined overset grid of size $h=1/200$ between $x=0.59$ and $x=0.91$ (where the solution contains a discontinuity) on a baseline grid of size $h=1/100$ as shown in Figure \ref{fig:SodOversetGrid}. The computed solution for density obtained at $t=0.2$ using the $h=1/200$ overset grid on $h=1/100$ baseline grid at the $y=0.5$ line for $P^{1}$, $P^{2}$ and $P^{3}$ based DGM is compared and plotted against the exact solution in Figure \ref{fig:SodOversetSolution}. We also plot the solution difference ($|\rho_{withOverset}-\rho_{withOutOverset}|$) obtained for $P^{1}$, $P^{2}$ and $P^{3}$ based DGM in Figure \ref{fig:SodOversetError} between the solution obtained using a grid of size $h=1/200$ without any overset and using a refined overset grid of size $h=1/200$ between $x=0.59$ and $x=0.91$ (where the solution contains a discontinuity) on a baseline grid of size $h=1/100$. From looking at Figure \ref{fig:SodOversetError}, we can see the solution obtained with a refined overset grid of size $h=1/200$ between $x=0.59$ and $x=0.91$ on a baseline grid of size $h=1/100$ is as good as the solution obtained with single grid of size $h=1/200$ especially on the overset grid.
\\
\\

\begin{figure}[htbp]
\begin{center}
\includegraphics[scale=0.2]{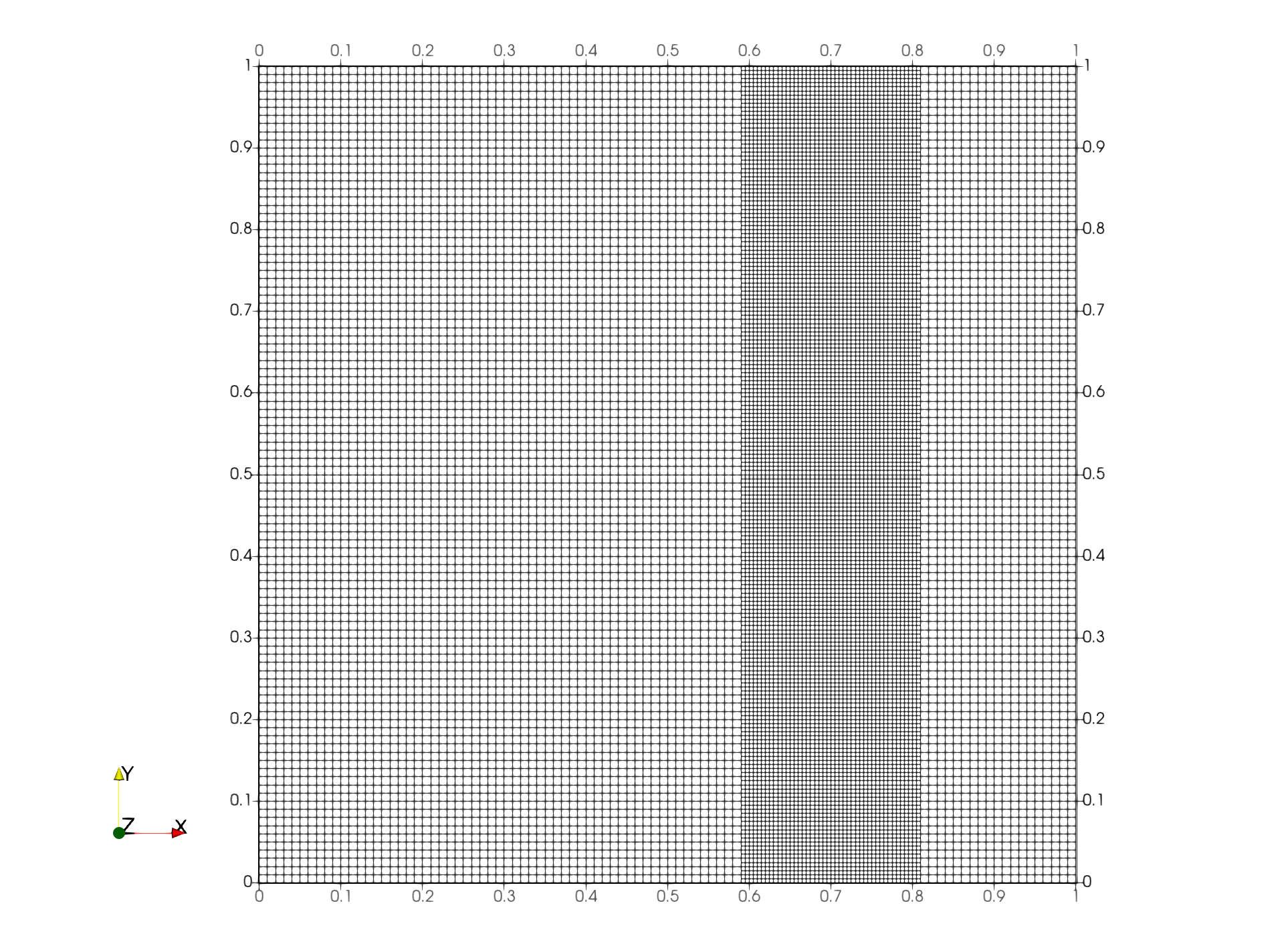}
\caption{Overset Grid used for the solution of Lax problem in the domain $[0,1]\times[0,1]$ with a refined overset grid of size $h=1/200$ between $x=0.59$ and $x=0.81$ (where the solution contains a discontinuity) on a baseline grid of size $h=1/100$.}
\label{fig:LaxOversetGrid}
\end{center}
\end{figure}

\begin{figure}[htbp]
\begin{center}
\includegraphics[scale=0.7]{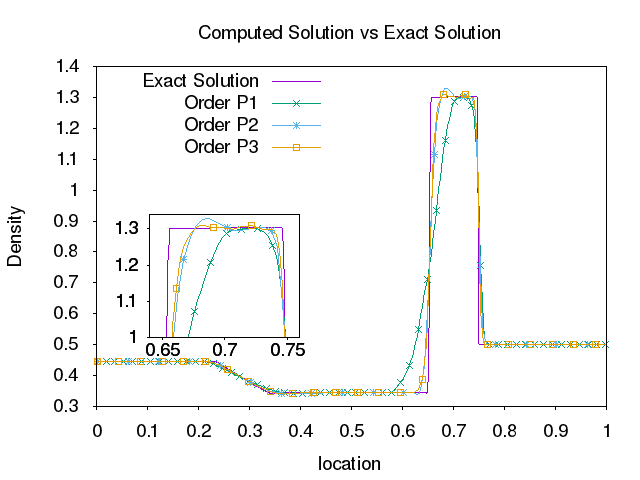}
\caption{Comparison of density solutions on $y=0.5$ line of Lax problem at $t=0.1$ in the domain $[0,1]\times[0,1]$ with a refined overset grid of size $h=1/200$ between $x=0.59$ and $x=0.81$ (where the solution contains a discontinuity) on a baseline grid of size $h=1/100$ obtained with the $P^{1}$, $P^{2}$ and $P^{3}$ based DGM using the proposed data communication scheme and the exact solution. Figure also includes a zoomed in portion of the solution for better comparison}
\label{fig:LaxOversetSolution}
\end{center}
\end{figure}

\begin{figure}[htbp]
\begin{center}
\includegraphics[scale=0.7]{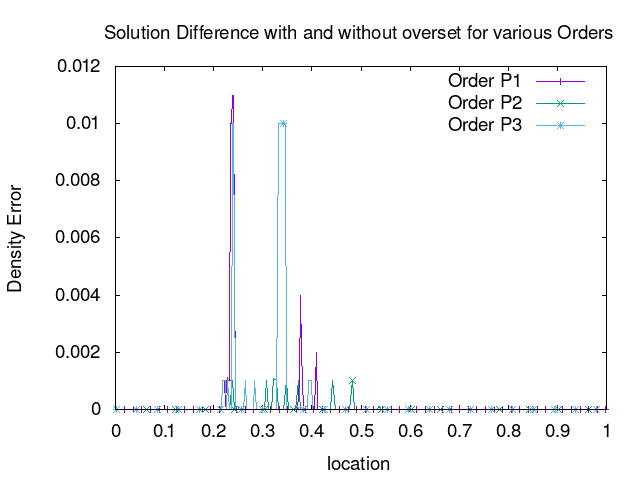}
\caption{Solution difference for density ($|\rho_{withOverset}-\rho_{withOutOverset}|$) obtained on $y=0.5$ line of Lax problem at $t=0.1$ in the domain $[0,1]\times[0,1]$ for $P^{1}$, $P^{2}$ and $P^{3}$ based DGM between a refined overset grid of size $h=1/200$ between $x=0.59$ and $x=0.81$ (where the solution contains a discontinuity) on a baseline grid of size $h=1/100$ and a single grid of size $h=1/200$.}
\label{fig:LaxOversetError}
\end{center}
\end{figure}

\noindent \textbf{Example 4:} We solve the Lax problem as proposed in \cite{lax1} in the two-dimensional domain. We solve the 2D Euler equations in the domain $[0,1] \times [0,1]$ with the initial conditions given as $(\rho,u,v,p)=(0.445,0.698,0.0,3.528)$ for $x < 0.5$ and $(\rho,u,v,p)=(0.5,0.0,0.0,0.571)$ otherwise. Non reflecting boundary condition is applied at $x=0$ and $x=1$ and periodic boundary conditions are applied at the other two boundaries. We use a refined overset grid of size $h=1/200$ between $x=0.59$ and $x=0.81$ (where the solution contains a discontinuity) on a baseline grid of size $h=1/100$ as shown in Figure \ref{fig:LaxOversetGrid}. The computed solution for density obtained at $t=0.1$ using the $h=1/200$ overset grid on $h=1/100$ baseline grid at the $y=0.5$ line for $P^{1}$, $P^{2}$ and $P^{3}$ based DGM is compared and plotted against the exact solution in Figure \ref{fig:LaxOversetSolution}. We also plot the solution difference ($|\rho_{withOverset}-\rho_{withOutOverset}|$) obtained for $P^{1}$, $P^{2}$ and $P^{3}$ based DGM in Figure \ref{fig:LaxOversetError} between the solution obtained using a grid of size $h=1/200$ without any overset and using a refined overset grid of size $h=1/200$ between $x=0.59$ and $x=0.81$ (where the solution contains a discontinuity) on a baseline grid of size $h=1/100$. From looking at Figure \ref{fig:LaxOversetError}, we can see the solution obtained with a refined overset grid of size $h=1/200$ between $x=0.59$ and $x=0.81$ on a baseline grid of size $h=1/100$ is as good as the solution obtained with single grid of size $h=1/200$ especially on the overset grid.
\\
\\

\noindent \textbf{Example 5:} As another test case, we look at the 2D Riemann problem of gas dynamics which is one of the most extensively studied problem which also contains a lot of intricate flow structures. We solve the two-dimensional Euler equations \eqref{2dEulerEquations} in the domain $[0,1]\times [0,1]$ for 2D Riemann problem configurations (10), (12) and (16) as given by the nomenclature in \cite{ll}. The initial conditions for configurations (10), (12) and (16) are given respectively as

\begin{flalign}\label{2DRiemannInitialConfig10}
 (\rho,u,v,p)(x,y,0) = \begin{cases}
                            (1,0,0.4297,1) \quad \text{if $x \geq 0.5$ and $y \geq 0.5$} \\
                            (0.5,0,0.6076,1) \quad \text{if $x<0.5$ and $y \geq 0.5$} \\
                            (0.2281,0,-0.6076,0.3333) \quad \text{if $x<0.5$ and $y<0.5$} \\
                            (0.4562,0,-0.4297,0.3333) \quad \text{otherwise}
                        \end{cases}
\end{flalign}

\begin{flalign}\label{2DRiemannInitialConfig12}
 (\rho,u,v,p)(x,y,0) = \begin{cases}
                            (0.5313,0,0,0.4) \quad \text{if $x \geq 0.5$ and $y \geq 0.5$} \\
                            (1,0.7276,0,1) \quad \text{if $x<0.5$ and $y \geq 0.5$} \\
                            (0.8,0,0,1) \quad \text{if $x<0.5$ and $y<0.5$} \\
                            (1,0,0.7276,1) \quad \text{otherwise}
                        \end{cases}
\end{flalign}

\begin{flalign}\label{2DRiemannInitialConfig16}
 (\rho,u,v,p)(x,y,0) = \begin{cases}
                            (0.5313,0.1,0.1,0.4) \quad \text{if $x \geq 0.5$ and $y \geq 0.5$} \\
                            (1.0222,-0.6179,0.1,1) \quad \text{if $x<0.5$ and $y \geq 0.5$} \\
                            (0.8,0.1,0.1,1) \quad \text{if $x<0.5$ and $y<0.5$} \\
                            (1,0.1,0.8276,1) \quad \text{otherwise}
                        \end{cases}
\end{flalign}

\begin{figure}[htbp]
\begin{center}
\includegraphics[scale=0.2]{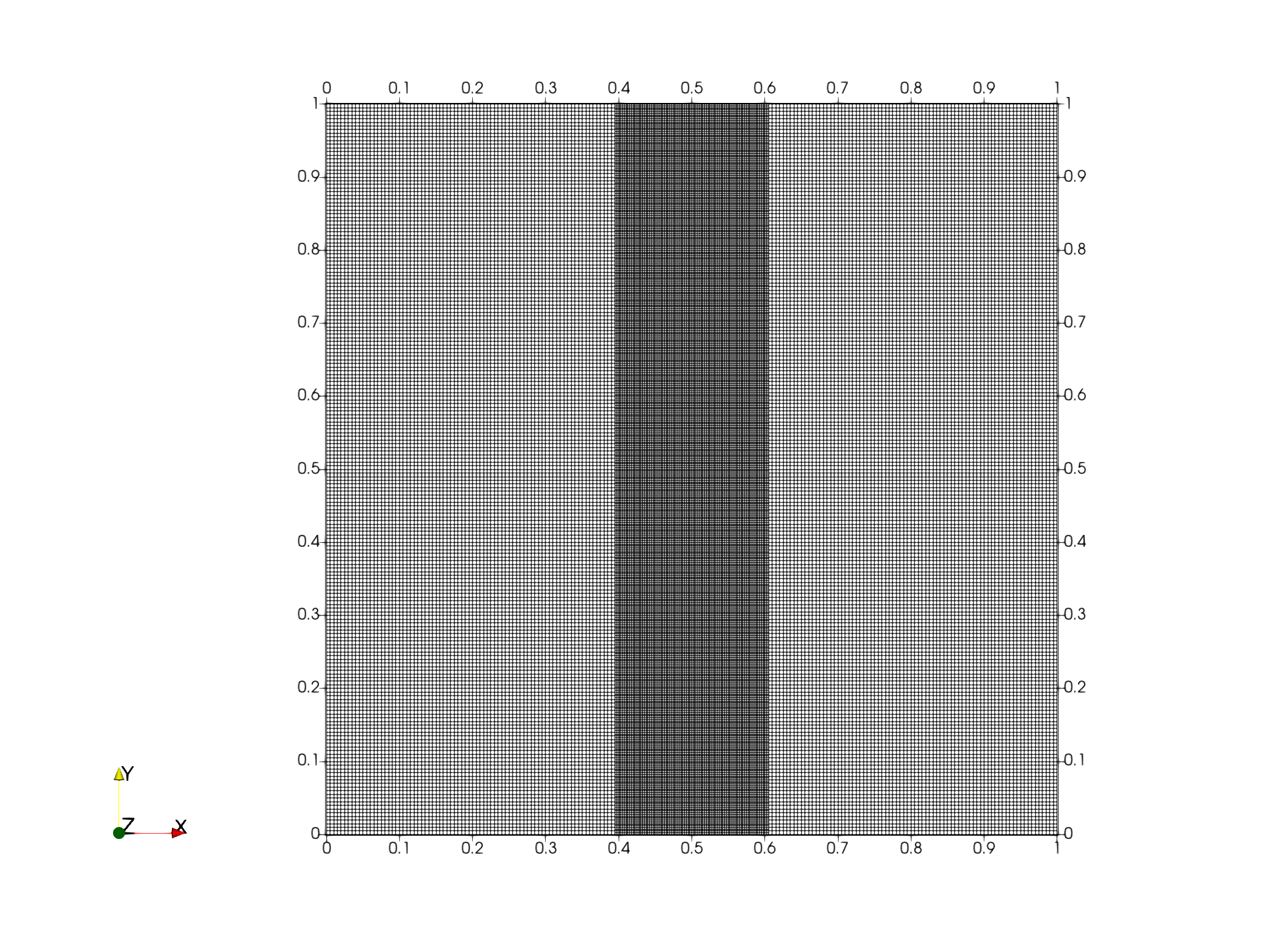}
\caption{Overset Grid used for the solution of 2D Riemann problem configurations (10) and (12) in the domain $[0,1]\times[0,1]$ with a refined overset grid of size $h=1/400$ between $x=0.395$ and $x=0.605$ (which is our region of interest) on a baseline grid of size $h=1/200$.}
\label{fig:2DRiemannConfig1012OversetGrid}
\end{center}
\end{figure}

\begin{figure}[htbp]
\begin{center}
\includegraphics[scale=0.2]{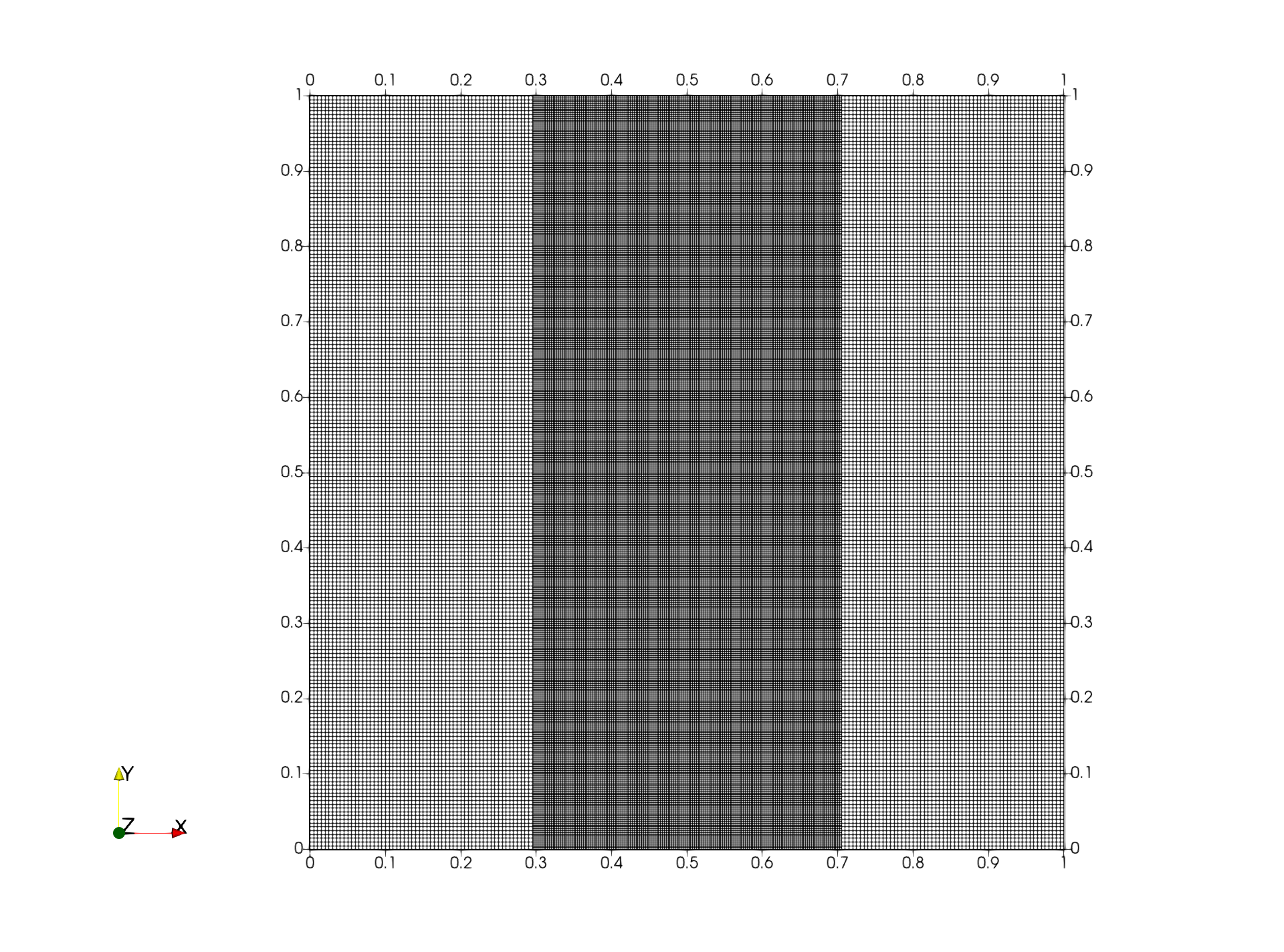}
\caption{Overset Grid used for the solution of 2D Riemann problem configuration (16) in the domain $[0,1]\times[0,1]$ with a refined overset grid of size $h=1/400$ between $x=0.295$ and $x=0.705$ (which is our region of interest) on a baseline grid of size $h=1/200$.}
\label{fig:2DRiemannConfig16OversetGrid}
\end{center}
\end{figure}

\begin{figure}[htbp]
\begin{center}
\includegraphics[width=0.95\textwidth]{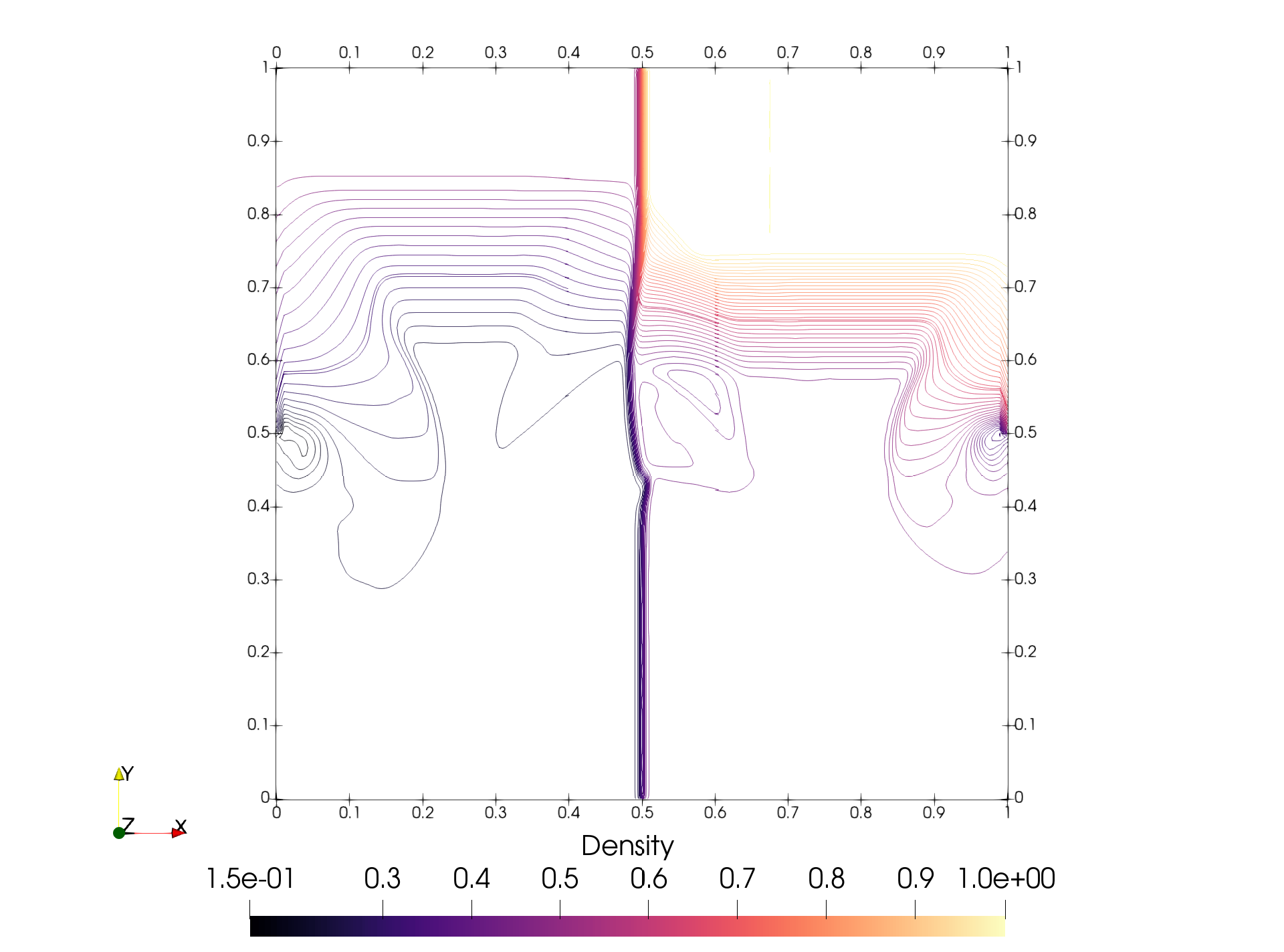}
\caption{50 equally spaced density contours for solution at $t=0.15$ for 2D Riemann problem configuration 10 with a refined overset grid of size $h=1/400$ between $x=0.395$ and $x=0.605$ (which is our region of interest) on a baseline grid of size $h=1/200$ obtained with the $P^{3}$ based DGM using the proposed data communication scheme}
\label{fig:2DRiemannConfig10OversetSolution}
\end{center}
\end{figure}

\begin{figure}[htbp]
\begin{center}
\includegraphics[width=0.95\textwidth]{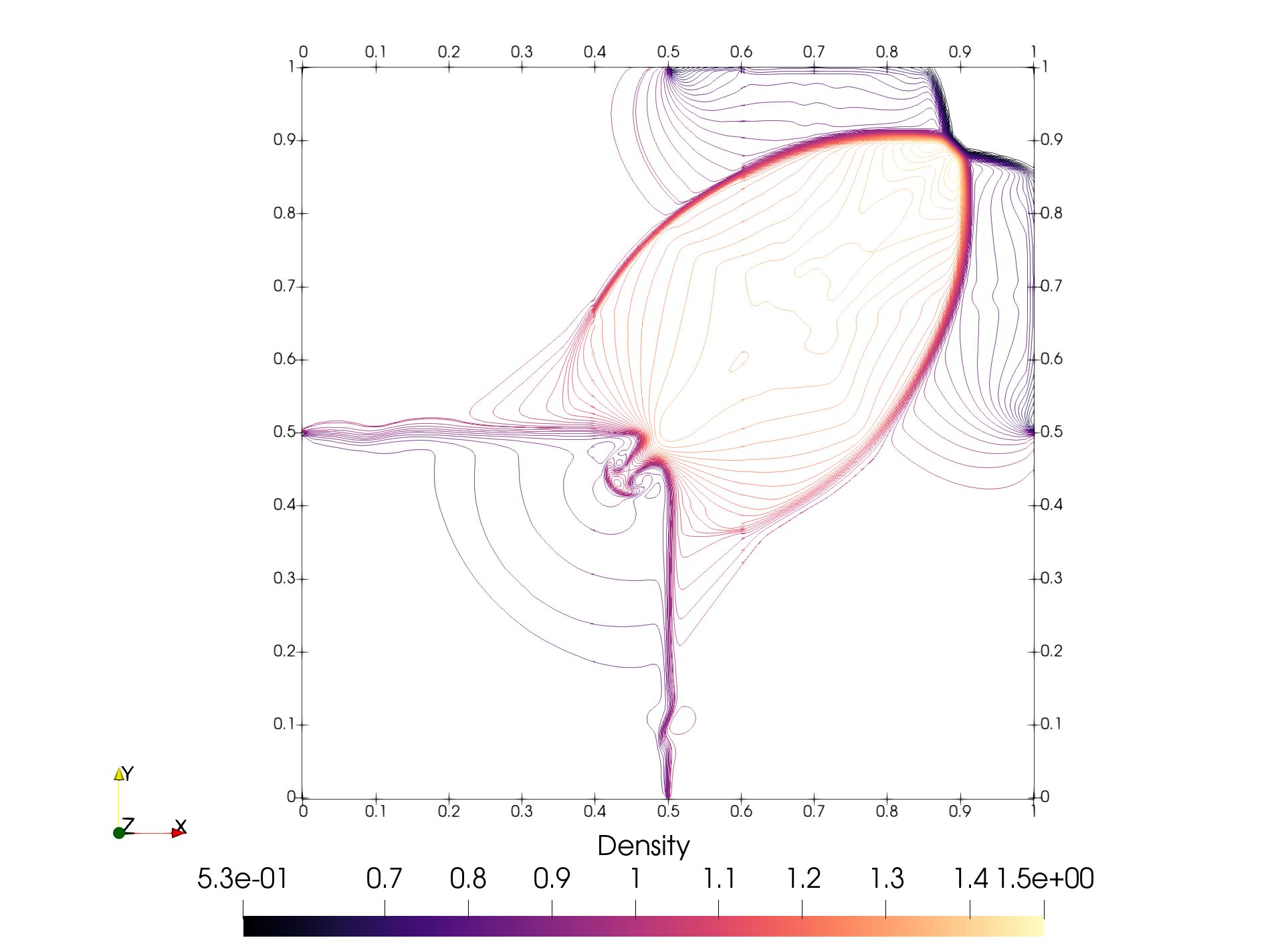}
\caption{50 equally spaced density contours for solution at $t=0.25$ for 2D Riemann problem configuration 12 with a refined overset grid of size $h=1/400$ between $x=0.395$ and $x=0.605$ (which is our region of interest) on a baseline grid of size $h=1/200$ obtained with the $P^{3}$ based DGM using the proposed data communication scheme}
\label{fig:2DRiemannConfig12OversetSolution}
\end{center}
\end{figure}

\begin{figure}[htbp]
\begin{center}
\includegraphics[width=0.95\textwidth]{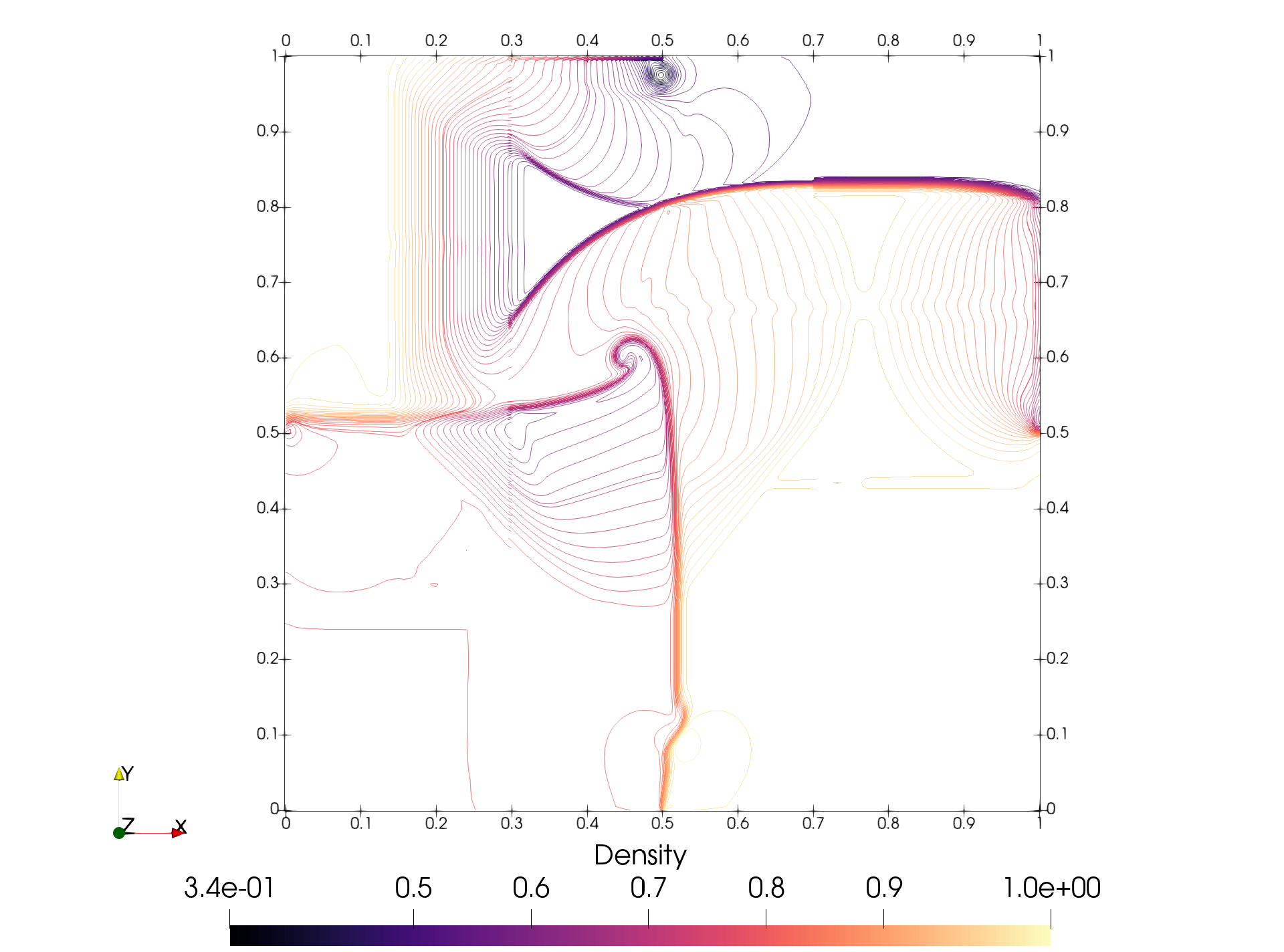}
\caption{50 equally spaced density contours for solution at $t=0.2$ for 2D Riemann problem configuration 16 with a refined overset grid of size $h=1/400$ between $x=0.395$ and $x=0.605$ (which is our region of interest) on a baseline grid of size $h=1/200$ obtained with the $P^{3}$ based DGM using the proposed data communication scheme}
\label{fig:2DRiemannConfig16OversetSolution}
\end{center}
\end{figure}

For configurations (10) and (12), we use a refined overset grid of size $h=1/400$ between $x=0.395$ and $x=0.605$ (which is our region of interest) on a baseline grid of size $h=1/400$ as shown in Figure \ref{fig:2DRiemannConfig1012OversetGrid}. For configuration (16), we use a refined overset grid of size $h=1/400$ between $x=0.295$ and $x=0.705$ (which is our region of interest) on a baseline grid of size $h=1/400$ as shown in Figure \ref{fig:2DRiemannConfig16OversetGrid}. We compute the solution upto time $t=0.15$ for configuration (10), $t=0.25$ for configuration (12) and till $t=0.2$ for configuration (16). Configuration (10) contains contact discontinuities and rarefaction waves initially. Configuration (12) and (16) contain both shocks and rarefaction waves along a contact discontinuities. We have selected the overset mesh for configurations (12) and (16) such that a shock passes through the overset grid. This will demonstrate that our data communication method also works as a limiter. We solved all three configurations for $\mathbf{P}^{1}$, $\mathbf{P}^{2}$ and $\mathbf{P}^{3}$ based DGM. The density contours for the solution obtained using the our procedure for $\mathbf{P}^{3}$ based DGM are shown in Figures \ref{fig:2DRiemannConfig10OversetSolution}, \ref{fig:2DRiemannConfig12OversetSolution} and \ref{fig:2DRiemannConfig16OversetSolution} respectively for configurations (10), (12) and (16). We also calculated the $L_{2}$ error for density for each of the solutions obtained for $\mathbf{P}^{1}$, $\mathbf{P}^{2}$ and $\mathbf{P}^{3}$ based DGM for all three configurations by taking a solution obtained on a single grid of size $h=1/400$ using $\mathbf{P}^{4}$ based DGM as the exact solution and this error is tabulated in Table \ref{table:3}. From the solution obtained, we can see that data communication scheme works quite well even when a shock passes through the artificial boundary of the overset grid (configurations (12) and (16)).

\begin{table}[htbp]
\large
\centering
\begin{tabular}{|c|c|c|c|}
\hline
\multicolumn{4}{|c|}{$L_{2}$ error for the three different Riemann problem configurations} \\
\hline
  & configuration (10) & configuration (12) & configuration (16) \\
\hline
$\mathbf{P}^{1}$ based DGM & 5.16E-05  & 7.49E-04 & 2.42E-04 \\
\hline
$\mathbf{P}^{2}$ based DGM & 6.63E-07  & 5.34E-06 & 9.94E-07\\
\hline
$\mathbf{P}^{3}$ based DGM & 8.83E-10  & 9.27E-09 & 2.76E-09\\
\hline
\end{tabular}
\caption{$L^{2}$ error for density obtained for 2D Riemann problem configurations (10), (12) and (16) using $\mathbf{P}^{1}$, $\mathbf{P}^{2}$ and $\mathbf{P}^{3}$ based DGM on the overset grids shown in Figures \ref{fig:2DRiemannConfig1012OversetGrid} and \ref{fig:2DRiemannConfig16OversetGrid} by using a solution obtained on a single grid of size $h=1/400$ using $\mathbf{P}^{4}$ based DGM as the exact solution.}
\label{table:3}
\end{table}

\section{Conclusion:}\label{sec:conc}

\noindent We have developed a new scheme for data communication using subcells and WENO reconstruction for two-dimensional problems using overset grids. We use element based data communication approach between overset grids and reconstruct the degrees of freedom in cells near the overset interface using WENO reconstruction. This is done by dividing the immediate neighbors into subcells as proposed in \cite{srspkmr1} and also by constructing a ghost cell near the artificial boundary of the element. This procedure has the added advantage that it also works as a limiter if a discontinuity passes through the overset interface. We can use this procedure for data communication between overset grids with any other higher order method which uses cells for their solution. We have demonstrated the scheme using discontinuous Galerkin method. We have provided accuracy tests to show that this procedure maintains the order of accuracy of the scheme on the overset grids. We have also provided results with solutions containing shocks to demonstrate the limiter aspect of this scheme.


\bibliographystyle{ieeetr}
\bibliography{references}

\end{document}

%% file: OversetExpl2New2.pdf_t
\begin{picture}(0,0)%
\includegraphics{OversetExpl2New2.pdf}%
\end{picture}%
\setlength{\unitlength}{4144sp}%
\begingroup\makeatletter\ifx\SetFigFont\undefined%
\gdef\SetFigFont#1#2#3#4#5{%
  \reset@font\fontsize{#1}{#2pt}%
  \fontfamily{#3}\fontseries{#4}\fontshape{#5}%
  \selectfont}%
\fi\endgroup%
\begin{picture}(6486,4998)(571,-4936)
\put(1216,-4921){\makebox(0,0)[lb]{\smash{{\SetFigFont{20}{24.0}{\familydefault}{\mddefault}{\updefault}{\color[rgb]{0,0,0}r}%
}}}}
\put(586,-4201){\makebox(0,0)[lb]{\smash{{\SetFigFont{20}{24.0}{\familydefault}{\mddefault}{\updefault}{\color[rgb]{0,0,0}s}%
}}}}
\put(631,-1096){\makebox(0,0)[lb]{\smash{{\SetFigFont{20}{24.0}{\familydefault}{\mddefault}{\updefault}{\color[rgb]{0,0,0}Grid 1}%
}}}}
\put(5041,-2806){\makebox(0,0)[lb]{\smash{{\SetFigFont{12}{14.4}{\familydefault}{\mddefault}{\updefault}{\color[rgb]{0,0,0}2}%
}}}}
\put(6166,-4786){\makebox(0,0)[lb]{\smash{{\SetFigFont{20}{24.0}{\familydefault}{\mddefault}{\updefault}{\color[rgb]{1,0,0}Grid 2}%
}}}}
\put(5446,-2491){\makebox(0,0)[lb]{\smash{{\SetFigFont{12}{14.4}{\familydefault}{\mddefault}{\updefault}{\color[rgb]{0,0,0}$\Omega_{k+1}$}%
}}}}
\put(4816,-2491){\makebox(0,0)[lb]{\smash{{\SetFigFont{12}{14.4}{\familydefault}{\mddefault}{\updefault}{\color[rgb]{0,0,0}$\Omega_{k}$}%
}}}}
\put(4636,-2806){\makebox(0,0)[lb]{\smash{{\SetFigFont{12}{14.4}{\familydefault}{\mddefault}{\updefault}{\color[rgb]{0,0,0}1}%
}}}}
\put(4636,-2311){\makebox(0,0)[lb]{\smash{{\SetFigFont{12}{14.4}{\familydefault}{\mddefault}{\updefault}{\color[rgb]{0,0,0}3}%
}}}}
\put(5086,-2311){\makebox(0,0)[lb]{\smash{{\SetFigFont{12}{14.4}{\familydefault}{\mddefault}{\updefault}{\color[rgb]{0,0,0}4}%
}}}}
\put(3646,-2491){\makebox(0,0)[lb]{\smash{{\SetFigFont{12}{14.4}{\familydefault}{\mddefault}{\updefault}{\color[rgb]{0,0,0}$\Omega_{k-1}$}%
}}}}
\put(4636,-3616){\makebox(0,0)[lb]{\smash{{\SetFigFont{12}{14.4}{\familydefault}{\mddefault}{\updefault}{\color[rgb]{0,0,0}$\Omega_{k-P}$}%
}}}}
\put(4636,-1771){\makebox(0,0)[lb]{\smash{{\SetFigFont{12}{14.4}{\familydefault}{\mddefault}{\updefault}{\color[rgb]{0,0,0}$\Omega_{k+P}$}%
}}}}
\end{picture}%

%% file: OversetExpl3New.pdf_t
\begin{picture}(0,0)%
\includegraphics{OversetExpl3New.pdf}%
\end{picture}%
\setlength{\unitlength}{4144sp}%
\begingroup\makeatletter\ifx\SetFigFont\undefined%
\gdef\SetFigFont#1#2#3#4#5{%
  \reset@font\fontsize{#1}{#2pt}%
  \fontfamily{#3}\fontseries{#4}\fontshape{#5}%
  \selectfont}%
\fi\endgroup%
\begin{picture}(6306,5223)(571,-4936)
\put(1216,-4921){\makebox(0,0)[lb]{\smash{{\SetFigFont{20}{24.0}{\familydefault}{\mddefault}{\updefault}{\color[rgb]{0,0,0}r}%
}}}}
\put(586,-4201){\makebox(0,0)[lb]{\smash{{\SetFigFont{20}{24.0}{\familydefault}{\mddefault}{\updefault}{\color[rgb]{0,0,0}s}%
}}}}
\put(631,-1096){\makebox(0,0)[lb]{\smash{{\SetFigFont{20}{24.0}{\familydefault}{\mddefault}{\updefault}{\color[rgb]{0,0,0}Grid 1}%
}}}}
\put(5986,-4561){\makebox(0,0)[lb]{\smash{{\SetFigFont{20}{24.0}{\familydefault}{\mddefault}{\updefault}{\color[rgb]{1,0,0}Grid 2}%
}}}}
\put(4501,-2311){\makebox(0,0)[lb]{\smash{{\SetFigFont{12}{14.4}{\familydefault}{\mddefault}{\updefault}{\color[rgb]{0,0,0}$\Omega_{k}$}%
}}}}
\end{picture}%